\documentclass{amsart}

\usepackage[curve]{xy}
\usepackage{amssymb}

\numberwithin{equation}{section}

\newtheorem{Thm}{Theorem}[section]

\theoremstyle{remark}
\newtheorem{Rem}[Thm]{Remark}

\theoremstyle{definition}

\begin{document}

\title[Punctured disks and the quantization of Poisson manifolds]
{Associative algebras, punctured disks and the quantization of Poisson 
manifolds}

\author{Domenico~Fiorenza}
\address{Dipartimento di Matematica ``G. Castelnuovo'',  
Universit\`a di Roma ``La Sapienza'', Piazzale Aldo Moro, 2 --- 
I-00185 Roma --- Italy}
\email{fiorenza@mat.uniroma1.it}

\author{Riccardo~Longoni}
\address{Dipartimento di Matematica ``G. Castelnuovo'',  
Universit\`a di Roma ``La Sapienza'', Piazzale Aldo Moro, 2 --- 
I-00185 Roma --- Italy}
\email{longoni@mat.uniroma1.it}

\begin{abstract}
The aim of the note is to provide an introduction to the algebraic, 
geometric and quantum field theoretic ideas that lie behind the 
Kontsevich-Cattaneo-Felder formula for the quantization of Poisson 
structures. We show how the quantization formula itself naturally 
arises when one imposes the following two requirements to a Feynman integral:
on the one side it has to reproduce the given Poisson structure as the first 
order term of its perturbative expansion; on the other side its 
three-point functions should describe an associative algebra. 
It is further shown how the Magri-Koszul brackets on 1-forms naturally fits 
into the theory of the Poisson sigma-model.
\end{abstract}

\maketitle

\section{Deformation quantization as a Feynman diagrams expansion}
\label{Fiorenza:s1}

A Poisson manifold is a differentiable manifold $M$ 
endowed with a bi-vector $\alpha\in\Gamma(M;TM\wedge TM)$ such that 
$[\alpha,\alpha]=0$, 
where $[,]$ is the Schouten-Nijenhuis bracket (see e.g. \cite{Fiorenza:V}).  
The bi-vector $\alpha$ defines a Poisson algebra structure on the space of 
smooth functions on $M$ by
\begin{displaymath}
\{f,g\}:=\langle\alpha|\mathrm d f\wedge \mathrm d g\rangle
\end{displaymath}
The problem of deformation quantization of the given Poisson structure 
is that of finding an associative $\star$-product on $C^\infty(M)[[\hbar]]$ 
deforming the usual pointwise product on $C^\infty(M)$ and having the 
Poisson bracket as the first order term in $\hbar$:
\begin{equation}\label{Fiorenza:deformation}
(f\star g)(x)=f(x)g(x)+\frac{{\tt i}\hbar}2 \{f,g\}(x)+O(\hbar^2),
\end{equation}
or, more generally,
\begin{equation}\label{Fiorenza:deformation2}
(f\star g)(x)=f(x)g(x)+\frac{{\tt
i}\hbar}2\biggl(\{f,g\}+B(f,g)\biggr)(x)+O(\hbar^2),
\end{equation}
where $B$ is a symmetric bi-differential operator. This problem has been 
solved by M.~Kontsevich \cite{Fiorenza:Kont}, and his solution was then 
interpreted in the language of quantum field theories by A.~Cattaneo 
and G.~Felder \cite{Fiorenza:CF1}. These notes are an attempt to explain 
why the Cattaneo-Felder model naturally arises when one tries to look
at \eqref{Fiorenza:deformation} as the perturbative expansion of 
a Feynman integral:
\begin{displaymath}
f\star g=\quad
\begin{xy}
,(0,0)*{\bullet},(0,-3)*{f},(5,0)*{\bullet},(5,-3)*{g}
\end{xy}\quad+\quad\frac{{\tt i\hbar}}{2}
\begin{xy}
,(0,7);(-5,-8)*{\bullet}**\dir{-} ?>*\dir{>} ,(-5,-11)*{f}
,(0,7);(5,-8)*{\bullet}**\dir{-} ?>*\dir{>} ,(5,-11)*{g}
,(0,9)*{\alpha}
\end{xy}\quad+O(\hbar^2)
\end{displaymath}
We see from this formula that there are two types of vertices, namely 
the ones labelled by the functions $f,g$ and the ones labelled by the 
bi-vector $\alpha$, and that the propagator is
\begin{displaymath}
\begin{xy}
,(-5,5);(5,5)**\crv{(-4,-5)&(4,-5)} ?>*\dir{>}
\end{xy}= \mathrm d x^i\otimes \partial_i
\end{displaymath}
where $\partial_i$ is a shorthand notation for $\partial/\partial x^i$.
By the above description, we see that our fields are tangent and 
cotangent vectors at $x$; moreover, in order to look at $\alpha$ as 
to a function of the fields, we have to consider the cotangent vectors 
as odd fields, i.e., the coordinates $\eta_i$ of a cotangent vector $\eta$ 
are anticommuting variables. Therefore, the natural choice for the space 
of fields is $T_xM\oplus\Pi T^*_xM$, endowed with 
the natural pairing $\langle \partial_i|{\mathrm d}x^j\rangle=\delta_i^j$.

The functions $f$ and $g$ and the Poisson bi-vector 
$\alpha$ can be seen as functions on the space 
of fields, by using the Taylor expansions:
\begin{align*}
f(\xi,\eta)&:=
f(x+\xi)=f(x)+\partial_if(x)\xi^i+\frac12\partial_i
\partial_jf(x)\xi^i\xi^j+\cdots\\
g(\xi,\eta)&:=
g(x+\xi)=g(x)+\partial_ig(x)\xi^i+\frac12\partial_i
\partial_jg(x)\xi^i\xi^j+\cdots\\
\alpha(\xi,\eta)&:=
\langle\alpha(x+\xi)|\eta\wedge\eta\rangle=\alpha^{ij}(x)\eta_i\eta_j +
\partial_k \alpha^{ij}(x)\eta_i\eta_j\xi^k+\cdots
\end{align*}
where $\xi\in T_xM$ and $\eta\in \Pi T^*_xM$. Now consider 
\begin{equation}\label{Fiorenza:perturbative}
\frac{\displaystyle{
\int_{T_xM\oplus\Pi T^*_xM} \mathrm d\xi \mathrm d
\eta\, f(x+\xi)\,g(x+\xi)\,
e^{\frac{{\tt i}}{\hbar}S(\xi,\eta)}
}}{\displaystyle{
\int_{T_xM\oplus\Pi T^*_xM} \mathrm d\xi \mathrm d
\eta\,
e^{\frac{{\tt i}}{\hbar}\langle\xi|\eta\rangle}
}}
\end{equation}
where the action is
\begin{displaymath}
S(\xi,\eta)=S_{\rm free}(\xi,\eta)+S_{\rm int}(\xi,\eta):=
\langle \xi|\eta\rangle+\langle 
\alpha(x+\xi)|\eta\wedge\eta\rangle.
\end{displaymath}
By the usual Feynman rules, the perturbative expansion 
of~\eqref{Fiorenza:perturbative} is
\begin{displaymath}
\begin{xy}
,(0,0)*{\bullet},(0,-3)*{f},(5,0)*{\bullet},(5,-3)*{g}
\end{xy}\quad+\quad\frac{{\tt i\hbar}}{2}\left(
\begin{xy}
,(0,7);(-5,-8)*{\bullet}**\dir{-} ?>*\dir{>} ,(-5,-11)*{f}
,(0,7);(5,-8)*{\bullet}**\dir{-} ?>*\dir{>} ,(5,-11)*{g}
,(0,9)*{\alpha}  
\end{xy}\quad+
\begin{xy}
,(-5,-8)*{\bullet},(-5,-11)*{f}
,(0,-1);(0,-8)*{\bullet}**\dir{-} ?>*\dir{>} ,(0,-11)*{g}
,(0,-1);(0,-1)**\crv{(-8,5)&(0,12)&(8,5)} ?>*\dir{>}
,(2,-2)*{\alpha}  
\end{xy}
\quad+
\begin{xy}
,(5,-8)*{\bullet},(5,-11)*{g}
,(0,-1);(0,-8)*{\bullet}**\dir{-} ?>*\dir{>} ,(0,-11)*{f}
,(0,-1);(0,-1)**\crv{(-8,5)&(0,12)&(8,5)} ?>*\dir{>}
,(2,-2)*{\alpha}  
\end{xy}
\right)\quad+O(\hbar^2)
\end{displaymath}
which is of the form \eqref{Fiorenza:deformation2}.
Note that, if $\alpha$ is constant as a function of $x\in M$, then the 
perturbative expansion of \eqref{Fiorenza:perturbative} is
\begin{displaymath}
(f\star g)(x)=\sum_{n=0}^\infty \frac{1}{n!}\left(\frac{{\tt i}\hbar}{2}
\right)^n \alpha^{i_1j_1}\cdots\alpha^{i_nj_n}\partial_{i_1}\cdots
\partial_{i_n}f(x)\, \partial_{j_1}\cdots
\partial_{j_n}g(x)
\end{displaymath}
which is precisely the Moyal $\star$-product formula. However, for general
$\alpha$, formula \eqref{Fiorenza:perturbative} does not yield an 
associative $\star$-product. A way to remedy this is to consider a 
topological space whose geometry describes the structure of associative 
algebras, and pull back our integral onto this space.

\section{Punctured disks and associative algebras}
Let $D$ be the unit complex disk, and let $B_{n}$ be the moduli space of
$(n+1)$ points on the boundary of $D$, for $n\ge 2$. The disk $D$ is 
identified with the complex upper half plane and its boundary with 
$\mathbb R\cup \{\infty\}$. Since the group of the biholomorphisms 
acts $3$-transitively on the set of boundary points on $D$, 
we can fix three of them to 
be $0$, $1$ and $\infty$, and make all the others lie in the interval 
$(0,1)$. Therefore $B_{n}$ is just the open $(n-2)$-dimensional simplex 
$0<t_1<\cdots<t_{n-2}<1$. One can define a compactification $\overline B_n$ 
of $B_n$ by adding products of 
$B_{n'}$, $n'<n$; these new boundary components correspond to the 
collapsing of two or more points in the boundary. For instance, there 
are two boundary components in $\overline B_3$ corresponding to the 
degenerations as $t=t_1$ goes to $0$ or to $1$.

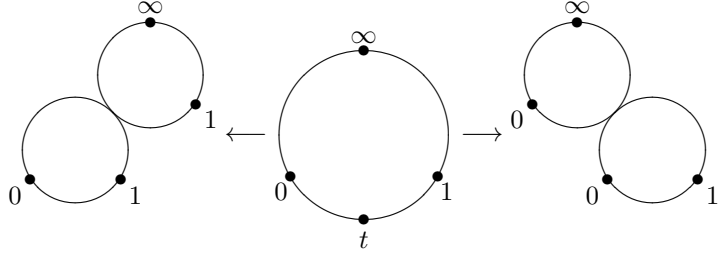
\begin{displaymath}
\begin{xy}
,(0,8)*\cir(8,8){},(-10,-2)*\cir(8,8){}
,(-16,-6)*{\bullet},(-4,-6)*{\bullet}
,(6,4)*{\bullet},(0,14.9)*{\bullet}
,(0,17)*{\infty},(-18,-8)*{0},(-2,-8)*{1}
,(8,2)*{1}
\end{xy}
\longleftarrow
\begin{xy}
*\cir(16,0){},(0,-11.3)*{\bullet},(-9.8,-5.65)*{\bullet}
,(9.8,-5.65)*{\bullet},(0,11.2)*{\bullet}
,(0,13)*{\infty},(0,-14)*{t},(-11,-7.5)*{0},(11,-7.5)*{1}
\end{xy}
\longrightarrow
\begin{xy}
,(0,8)*\cir(8,8){},(10,-2)*\cir(8,8){}
,(16,-6)*{\bullet},(4,-6)*{\bullet}
,(-6,4)*{\bullet},(0,14.9)*{\bullet}
,(0,17)*{\infty},(18,-8)*{1},(2,-8)*{0}
,(-8,2)*{0}
\end{xy}
\end{displaymath}

Now, we look at $B_2$ as to an operation $m_2$ with two inputs (the 
points $0$ and $1$) and one output (the point $\infty$). Note that
the two boundary components of $\overline B_3$ correspond to the
two ways of composing $m_2$ with itself, namely $m_2(m_2\otimes id)$
and $m_2(id\otimes m_2)$. So, if we find a continuous family of 
operations $m_3(t)$, $t\in(0,1)$, with three inputs and one output, 
which extends to the compactification $\overline B_3$ (in a way compatible 
with the product structure of the boundary), then the associativity 
of $m_2$ is equivalent to $m_3(0) =m_3(1)$. 
If moreover $m_3(t)$ is differentiable, this is equivalent to  
\begin{displaymath}
m_2\,\,\,\text{associative}\quad\Leftrightarrow \quad
\int_0^1\mathrm d t \,\frac{\mathrm dm_3(t)}{\mathrm d t} = 0
\end{displaymath}

\begin{Rem}
In the language of operads, the above discussion corresponds to
the well-known fact that the chain complex $C_*(\overline B_n)$ is the 
operad governing $A_\infty$ algebras. In particular one says that 
$m_2$ is associative only up to the homotopy $m_3$. 
\end{Rem}

Now, we want to define $m_2$ and $m_3$ on the space of 
smooth functions on the Poisson manifold $M$, in such a way that $m_2$ is 
related to eq.~\eqref{Fiorenza:perturbative}. The most natural choice is
to consider the ``expectation value'' over the maps $X\colon D\to M$ 
of the product $f(X(0))\,g(X(1))\,h(X(\infty))$ w.r.t. some measure 
to be defined, and ``raise'' the indices, i.e., set $h$ to be the Dirac
delta function $\delta_x$. In other words we are looking for an 
operation $m_2$ of the form
\begin{equation}
\label{Fiorenza:m2}
m_2(f,g)(x) = \int\mathrm d\mu(X)\, f(X(0))\,g(X(1))\,
\delta_x(X(\infty)).
\end{equation}
As for $m_3=m_3(t)$, we set
\begin{displaymath}
m_3(f,g,h)(x) = \int\mathrm d\mu(X)\, f(X(0))\,g(X(t))\,h(X(1))\,
\delta_x(X(\infty))
\end{displaymath}
so that the associativity of $m_2$ becomes
\begin{equation}
\label{Fiorenza:wbw}
\int\mathrm d\mu(X)\int_0^1\mathrm dt \left(\, f(X(0))\,
\frac{\mathrm dg(X(t))}{\mathrm dt}\,h(X(1))\,\delta_x(X(\infty))
\right)=0.
\end{equation}

\section{The Poisson sigma-model}

In this Section we want to combine eq.~\eqref{Fiorenza:m2}, which defines 
an associative product, with eq.~\eqref{Fiorenza:perturbative}, which has 
the correct first term in its perturbative expansion. First, the
measure $\mathrm d\mu(X)$ in eq.~\eqref{Fiorenza:m2} should be of the form
$(1/C)\,\mathrm d\xi \mathrm d e^{\frac{{\tt i}}{\hbar}S(\xi,\eta)}$ as in 
eq.~\eqref{Fiorenza:perturbative}, where $C$ is a suitable normalization
constant. In order to accomplish this, a new 
field, denoted by $\eta$, has to be introduced: it has to be defined 
on the disk and take values in $\Pi T^*_xM$. Moreover, since the new 
action $S$ will be an integral over $D$, it is natural to take 
$\eta\in\Omega^1(D;X^*(\Pi T^*M))$. We are therefore led to consider the 
following object
\begin{equation}\label{Fiorenza:CF}
\frac{\displaystyle{\int \mathrm dX \mathrm d
\eta\, f(X(0))\,g(X(1))\delta_x(X(\infty))\,
e^{\frac{{\tt i}}{\hbar}S(X,\eta)}
}}{\displaystyle{
\int \mathrm dX \mathrm d
\eta\,
e^{\frac{{\tt i}}{\hbar}\int_D\langle\mathrm dX|\eta\rangle}
}}
\end{equation}
where $S(X,\eta) = \int_D\langle \mathrm d X|\eta\rangle+\frac12 \int_D\langle 
\alpha(X)|\eta\wedge\eta\rangle$.

Notice however that in eq.~\eqref{Fiorenza:perturbative}, we have a tangent
vector $\xi\in T_xM$, where $x$ is some point in $M$. Hence, what we 
should consider are infinitesimal variations of the map $X$ around the 
constant map $X\equiv x$. In other terms, in eq.~\eqref{Fiorenza:CF}
we have to replace $X$ with $x+\xi$ where $\xi\in\Omega^0(D;X^*(TM))$.

Since the map $X$ at the point $\infty$ is fixed to be equal to $x$ by the 
term $\delta_x(X(\infty))$, we have to impose the boundary condition 
$\xi(\infty)=0$; finally the $1$-form $\eta$ is required to vanish on 
tangent vectors to the boundary of the disk $D$.
The action now reads
\begin{displaymath}
S(\xi,\eta)=\int_D\langle\mathrm d\xi|\eta\rangle+\frac12 \int_D\langle 
\alpha(x+\xi)|\eta\wedge\eta\rangle,
\end{displaymath}
and we define 
\begin{equation}\label{Fiorenza:star}
(f\star g)(x):=
\frac{\displaystyle{
\int \mathrm d\xi \mathrm d
\eta\, f(x+\xi(0))\,g(x+\xi(1))\,
e^{\frac{{\tt i}}{\hbar}S(\xi,\eta)}
}}{\displaystyle{
\int \mathrm d\xi \mathrm d
\eta\,
e^{\frac{{\tt i}}{\hbar}\int_D\langle \mathrm d\xi|\eta\rangle}
}}.
\end{equation} 

In order to perform the perturbative expansion of \eqref{Fiorenza:star},
symmetries of the action have to be taken into account. A systematic way 
of doing this is via the superfield formalism, namely we consider the 
superdisk $D^{2|2}$ with even coordinates $u^1,u^2$ and 
Grassmann coordinates $\theta^1,\theta^2$ and set
\begin{align*}
\tilde \xi^i&=\xi^i+\eta^{+\,i}_\mu\theta^\mu+\frac12\beta^{+\,i}_{\mu\nu}
\theta^\mu\theta^\nu\\
\tilde \eta_i&=\beta_i+\eta_{i\,\mu}\theta^\mu+\frac12\xi^{+}_{i\,\mu\nu}
\theta^\mu\theta^\nu.
\end{align*}
The de~Rham differential now reads $\mathrm D=\theta_\mu\frac{\partial}
{\partial u^\mu}$ and the $\star$-product becomes
\begin{equation}\label{Fiorenza:star-tilde}
(f\star g)(x)=
\frac{\displaystyle{
\int_{\xi^+=\eta^+=\beta^+=0} \mathrm d\tilde\xi \mathrm d
\tilde\eta\, f(x+\tilde\xi(0))\,g(x+\tilde\xi(1))\,
e^{\frac{{\tt i}}{\hbar}S(\tilde\xi,\tilde\eta)}
}}{\displaystyle{
\int_{\xi^+=\eta^+=\beta^+=0} \mathrm d\tilde\xi \mathrm d
\tilde\eta\,
e^{\frac{{\tt i}}{\hbar}\int_{D^{2|2}} \langle \mathrm D
\tilde \xi|\tilde \eta\rangle} 
}},
\end{equation}
where the superaction is 
\begin{equation}\label{Fiorenza:action}
S(\tilde\xi,\tilde\eta):=\int_{D^{2|2}} \langle \mathrm D 
\tilde \xi|\tilde \eta\rangle+
\frac12 \int_{D^{2|2}} \langle \alpha(x+\tilde \xi)|\tilde\eta\wedge
\tilde\eta \rangle.
\end{equation}
Notice that besides of the original fields $\xi,\eta$ (and their 
``antifields'' $\xi^+,\eta^+$), a new field $\beta$ has appeared, which 
can be interpreted as an infinitesimal symmetry of the original 
action (see Remark~\ref{Fiorenza:rem1} below). 

The advantage of this reformulation of the Poisson sigma-model is that 
we can now apply the Batalin-Vilkovisky formalism and deform the subspace 
$\xi^+=\eta^+=\beta^+=0$ over which the integration is performed, in such 
a way that the perturbative expansion is well defined.

\section{Batalin-Vilkovisky formalism} 

We recall that for any vector space $V$, the space of functions on 
$V\oplus\Pi V^*$ is naturally endowed with a BV algebra structure 
\cite{Fiorenza:ASKZ, Fiorenza:Sch}. Using the standard terminology, we call 
{\em fields}\/ the coordinates $v^i$ on $V$ and {\em antifields}\/ 
the coordinates $v^+_i$ on $\Pi V^*$. The BV bracket between two 
functionals $f,g\colon V\oplus\Pi V^*\to \mathbb R$ is given by
\begin{displaymath}
(f,g):=\frac{\overleftarrow\partial f}{\partial v^i}
\frac{\overrightarrow\partial g}{\partial v^+_i} - 
\frac{\overleftarrow\partial f}{\partial v^+_i}
\frac{\overrightarrow\partial g}{\partial v^i}
\end{displaymath}
while the BV Laplacian is
\begin{displaymath}
\Delta f = 
\frac{\overrightarrow\partial}{\partial v^+_i}
\frac{\overleftarrow\partial}{\partial v^i} f
\end{displaymath}
The BV bracket and the BV Laplacian satisfy, together with the
pointwise product, the axioms of a BV algebra, namely 
\begin{align*}
&(f,g) = -(-1)^{(|f|-1)\,(|g|-1)}(g,f)\\
& (f,(g,h)) = ((f,g),h) + (-1)^{(|f|-1)(|g|-1)}(g,(f,h))=0\\
&(f,gh)=(f,g)h + (-1)^{(|f|-1)|g|} g(f,h)\\
&(f,g) = \Delta(fg) - \Delta(f)g + (-1)^{|f|} f\Delta(g) \\ 
&\Delta^2=0
\end{align*}
In particular a $\Delta$-cohomology is defined on the space of 
functional on the fields-antifields.

In our case 
\begin{align*}
(\xi,\eta,\beta)&\in V=
\Omega^0(D,X^*(TM))\oplus \Omega^1(D,X^*(\Pi T^*M))
\oplus \Omega^0(D,X^*(\Pi T^*M))\\
(\xi^+,\eta^+,\beta^+)&\in \Pi V^*=
\Omega^2(D,X^*(\Pi T^*M))\oplus \Omega^1(D,X^*(TM))
\oplus \Omega^2(D,X^*(TM)). 
\end{align*}
A ``total degree'' is then introduced by 
saying that a form on $D$ with values in $X^*(TM)$ has total 
degree zero, while a form with values in $X^*(\Pi T^*M)$ has total 
degree 1. Next, we define the ``ghost number'' $\mathit{gh}$ as the 
difference between the total degree and the degree $\mathit{deg}$ 
as a differential form on $D$. We summarize the degrees and ghost numbers 
of our fields and antifields in the following table:\\
\begin{tabular}[t]{|c|c|c|c|c|c}
\hline$\mathit{gh}\backslash\mathit{deg}$ & 0 & 1 & 2 \\\hline
-2\hspace{0.5cm} &         &           & $\beta^+$ \\\hline
-1\hspace{0.5cm} &         & $\eta^+$  & $\xi^+$   \\ \hline
0\hspace{0.5cm}  & $\xi$   & $\eta$    & \\ \hline
1\hspace{0.5cm}  & $\beta$ &           &    \\\hline         
\end{tabular}
\vskip 11pt

A main feature of the BV formalism is that the integral
of a $\Delta$-closed functional $\mathcal H$ performed over a 
Lagrangian submanifold $L$ in the space of fields-antifields, depends 
only on the homology class of $L$ and that the integral of a 
$\Delta$-exact functional is zero. Hence, integration defines a 
pairing between homology classes of Lagrangian submanifolds and 
$\Delta$-cohomology classes. An easy computation shows that a 
functional of the form $e^{\frac{\tt i}\hbar S}$ is $\Delta$-closed 
if and only if $S$ satisfies the ``quantum master equation''
\begin{equation}
\label{Fiorenza:qmaster}
(S,S)-2i\hbar\Delta(S)=0
\end{equation}
as indeed happens for the superaction~\eqref{Fiorenza:action} of 
the Poisson sigma-model \cite{Fiorenza:CF1} (see also 
Remark~\ref{Fiorenza:rem-qme} below). More generally, if 
the functional $\mathcal H$ is of the form $\mathcal O\, 
e^{\frac{\tt i}\hbar S}$ for some functional $\mathcal O$ and 
some $S$ satisfying eq.~\eqref{Fiorenza:qmaster}, we have that 
$\Delta(\mathcal O\, e^{\frac{\tt i}\hbar S})=0$ if and only if 
$\Omega(\mathcal O)=0$, where $\Omega (\mathcal O)
:=(S,\mathcal O) -i\hbar\Delta(\mathcal O)$. 
Equation~\eqref{Fiorenza:qmaster} immediately implies $\Omega^2=0$ and 
the relevant cohomology classes are called ``observables'' of the theory.
Since the ``expectation value'' $\langle\mathcal O\rangle
:=\int_L \mathcal O\, e^{\frac{\tt i}\hbar S}$ of an observable $\mathcal O$  
depends only on the homology class of $L$,
the perturbative expansion of the original path 
integral~\eqref{Fiorenza:star}, which corresponds to integrating over 
the Lagrangian submanifold $\xi^+=\eta^+=\beta^+=0$ (and which is actually 
ill-defined due to the symmetries), can be effectively computed by choosing 
an appropriate submanifold where the quadratic part of the action is 
non-degenerate (see \cite{Fiorenza:CF1} for details).

\begin{Rem}
For any point $u$ in the boundary of $D$, one has
\begin{equation}
\label{Fiorenza:bordo}
\Omega(\tilde \xi^i(u))=\Omega(\tilde \eta_j(u))=0.
\end{equation}
This gives a way to construct observables for the Poisson sigma-model from
a point $u\in \partial D$ and a smooth function $\varphi$ of
$\tilde\xi$
and $\tilde\eta$. Indeed, the functional $\mathcal O_{\varphi,\,u}(\tilde
\xi,\tilde\eta)
:=\varphi(\tilde\xi(u),\tilde\eta(u))$ is clearly $\Omega$-closed.
In particular, $f(x+\tilde\xi(0))$ and $g(x+\tilde\xi(1))$ from 
eq.~\eqref{Fiorenza:star-tilde} are observables.  
\end{Rem}

\begin{Rem}\label{Fiorenza:rem-qme}
Given a $p$-multivector field $\psi$, written in coordinates as 
$\psi(x)^{i_1,\ldots, i_p} \partial_{i_1}\wedge\cdots\wedge 
\partial_{i_p}$, we can consider
\begin{displaymath}
S_{\psi}(\tilde\xi,\tilde\eta):= 
\int_{D^{2|2}}\psi(x+\tilde\xi)^{i_1,\ldots, i_p}
\tilde\eta_{i_1} \cdots\tilde\eta_{i_p}
\end{displaymath}
Notice that with this notation the superaction \eqref{Fiorenza:action}
becomes $S(\tilde\xi,\tilde\eta) = S_\text{free}(\tilde\xi,\tilde\eta)
+ S_\alpha(\tilde\xi,\tilde\eta)$. An explicit calculation shows that
the map $\psi\mapsto S_\psi$ is a Lie algebra morphism
\begin{displaymath}
(S_{\psi_1},S_{\psi_2})=S_{[\psi_1,\psi_2]}
\end{displaymath}
where we have the BV bracket on the l.h.s. and the Schouten-Nijenhuis
bracket on the r.h.s. In particular, since the bi-vector $\alpha$ is 
Poisson, we have $(S_\alpha,S_\alpha)=0$. When the ``free'' part of
the superaction is taken into account, it is not difficult to show that
$(S_\text{free},S_\text{free})=0$ and $(S_\text{free}, S_\psi)=0$, 
which in turn imply the so-called ``master equation''
for~\eqref{Fiorenza:action} 
\begin{equation}\label{Fiorenza:me}
(S,S)=(S_\text{free}+S_\alpha, S_\text{free}+S_\alpha)=0.
\end{equation}
A consequence of this equality is that $\delta f:=(S,f)$ 
is a coboundary operator. Finally, notice that the quantum master 
equation~\eqref{Fiorenza:qmaster} descends immediately from the relations 
$\Delta(S_\text{free}) = \Delta(S_\alpha) = 0$.
\end{Rem}

\begin{Rem}\label{Fiorenza:rem1}
Using the operator $\delta$ defined above, we can rewrite 
equation~\eqref{Fiorenza:me} as
\begin{equation}\label{Fiorenza:master}
\delta S =0.
\end{equation}
On the other hand one explicitly computes
\begin{align}\label{Fiorenza:symm1}
\delta \tilde\xi^i &= \mathrm D \tilde\xi^i + 
\alpha^{ij}(x+\tilde\xi)\tilde\eta_j\\
\label{Fiorenza:symm2}
\delta \tilde\eta_i &= \mathrm D \tilde\eta_i +\frac12 \partial_i\alpha^{jk}
(x+\tilde\xi)\tilde\eta_j\tilde\eta_k.
\end{align}
The operator $\delta|_{\xi^+=\eta^+=\beta^+=0}$ can be seen as a 
vector field on the space of functionals of $(\xi,\eta)$ depending on 
the choice of $\beta$. We denote by $\delta_\beta$ this vector field. 
Now, equations~ (\ref{Fiorenza:master}--\ref{Fiorenza:symm2}) 
together imply that $\delta_\beta$ is an 
infinitesimal symmetry of the original action $S(\xi,\eta)$. Explicitly 
this symmetry reads
\begin{align*}
\delta_\beta\xi^i &= \alpha^{ij}(x+\xi)\beta_j\\
\delta_\beta\eta_i &= -\mathrm d\beta_i -\partial_i\alpha^{jk}(x+\xi)
\eta_j\beta_k.  
\end{align*}
\end{Rem}

\section{Ward identities}\label{Fiorenza:trick}
The equation $\int_L\Delta(\mathcal H)=0$ produces non-trivial identities 
(called ``Ward identities'') among the expectation values. For instance 
if $\phi(\tilde\xi,\tilde\eta)$ is a $\Delta$-closed functional, the 
following equality easily descends from the axioms of a BV algebra
\begin{equation}
\label{Fiorenza:ward}
0=\int_L \Delta\left(e^{\frac{\tt i}\hbar S}
\phi\right) =
\int_L e^{\frac{\tt i}\hbar S} \delta\phi. 
\end{equation} 
Now consider
\begin{displaymath} 
\phi = \int_0^1\mathrm dt\int \mathrm d\theta\,
f(x+\tilde\xi(0))\,g(x+\tilde\xi(t,\theta))\,h(x+\tilde\xi(1)).
\end{displaymath} 
An explicit computation using eq.~\eqref{Fiorenza:symm1} shows that
\begin{displaymath} 
\delta\phi = \int_0^1\mathrm dt\,
\left(f(x+\tilde\xi(0))\,\frac{\mathrm d g(x+\tilde\xi(t))}{\mathrm dt}
\,h(x+\tilde\xi(1))\right),
\end{displaymath} 
Therefore eq.~\eqref{Fiorenza:ward} has precisely the form of 
eq.~\eqref{Fiorenza:wbw} and the Ward identity for this choice of 
$\phi$ is the associativity equation
\begin{displaymath}  
0=\int_L \Delta\left(e^{\frac{\tt i}\hbar S}
\phi\right) = ((f\star g)\star h) -(f\star (g\star h))
\end{displaymath}

\section{The Magri-Koszul bracket}

If $\omega$ is a $1$-form on $M$ we can associate to it a
function on $T_xM\oplus\Pi T^*_xM$ by  
\begin{displaymath}
\omega(\xi,\eta)=\omega(\xi):=\langle\omega(x+\xi)|
\xi\rangle=\omega_i(x)\xi^i+\partial_j\omega_i(x)\xi^i\xi^j+\cdots,
\end{displaymath}
Similarly, to a vector field $\chi$ we can associate the
function
\begin{displaymath}
\chi(\xi,\eta):=\langle\chi(x+\xi)|
\eta\rangle=\chi^i(x)\eta_i+\partial_j\chi^i(x)\eta_i\xi^j+\cdots,
\end{displaymath}
The perturbative expansion of the integral
\begin{displaymath}
\frac{\displaystyle{
\int_{T_xM\oplus\Pi T^*_xM}\mathrm d \xi\mathrm d\eta\,
\omega_1\bigl(\xi\bigr)\,
\omega_2\bigl(\xi\bigr)\,\chi\bigl(\xi,\eta\bigr)\,
e^{\frac{\tt i}{\hbar}S(\xi,\eta)}
}}{\displaystyle{
\int_{T_xM\oplus\Pi T^*_xM}\mathrm d \xi\mathrm d\eta\,
e^{\frac{\tt i}{\hbar}\langle\xi|\eta\rangle}
}}
\end{displaymath}
is closely related to the Magri-Koszul bracket on $1$-forms
\cite{Fiorenza:M,Fiorenza:Kosz}. More precisely, if we apply the Poisson 
sigma-model techniques to this situation, the function 
$\omega_1\bigl(\xi\bigr)\, \omega_2\bigl(\xi\bigr)\,
\chi\bigl(\xi,\eta\bigr)$ is changed into 
$\omega_1\bigl(\tilde\xi(0)\bigr)\, 
\omega_2\bigl(\tilde\xi(1)\bigr)\,
\chi\bigl(\tilde\xi(\infty),\tilde\eta(\infty)\bigr)$. 
Since $\xi(\infty)=0$, we have $\chi(\tilde\xi(\infty),\tilde\eta(\infty))
=\chi^i(x)\tilde\eta_i(\infty)$. Therefore the perturbative expansion of the 
path integral:
\begin{equation}\label{Fiorenza:MKintegral}
\frac{\displaystyle{
\int_{\xi^+=\eta^+=\beta^+=0} \mathrm d\tilde\xi \mathrm d \tilde\eta\, 
\omega_1(\tilde\xi(0))\,\omega_2(\tilde\xi(1))\,
\chi(\tilde\xi(\infty),\tilde\eta(\infty))\,
e^{\frac{{\tt i}}{\hbar}S(\tilde\xi,\tilde\eta)}
}}{\displaystyle{
\int_{\xi^+=\eta^+=\beta^+=0} \mathrm d\tilde\xi \mathrm d \tilde\eta\,
e^{\frac{{\tt i}}{\hbar}
\int_{D^{2|2}} \langle \mathrm D
\tilde \xi|\tilde \eta\rangle}
}}
,
\end{equation}
will depend on $\chi(x)$ but not on its derivatives.
The first order expansion of the integral~\eqref{Fiorenza:MKintegral} is
\begin{math}
\frac{{\tt i}\hbar}{2}\,
\bigl\langle\omega_1\bullet\omega_2\bigr|\chi\bigr\rangle
+O(\hbar^2)
\end{math}
where
\begin{displaymath}
\bigl\langle\omega_1\bullet\omega_2\bigr|\chi\bigr\rangle:=\,\,
\begin{xy}
,(0,5);(-5,-4)*{\bullet}**\dir{-} ?>*\dir{>}
,(-5,-7)*{\omega_1} ,(0,5);(5,-4)*{\bullet}**\dir{-}
?>*\dir{>} 
,(5,-7)*{\omega_2}
,(0,7)*{\alpha}
,(-8,3);(-5,-4)*{\bullet}**\dir{-} ?>*\dir{>}
,(-8,5)*{\chi}
\end{xy}
\quad + \quad
\begin{xy}
,(0,5);(-5,-4)*{\bullet}**\dir{-} ?>*\dir{>}
,(-5,-7)*{\omega_1} ,(0,5);(5,-4)*{\bullet}**\dir{-}
?>*\dir{>} 
,(5,-7)*{\omega_2}
,(0,7)*{\alpha}
,(8,3);(5,-4)*{\bullet}**\dir{-} ?>*\dir{>}
,(8,5)*{\chi}
\end{xy}
\quad + \quad 
\begin{xy}
,(0,1);(-5,-4)*{\bullet}**\dir{-} ?>*\dir{>}
,(-5,-7)*{\omega_1} ,(0,1);(5,-4)*{\bullet}**\dir{-}
?>*\dir{>} 
,(5,-7)*{\omega_2}
,(0,7)*{\chi}
,(0,5);(0,1)**\dir{-} ?>*\dir{>}
,(2,1)*{\alpha}
\end{xy}
\quad + \quad 
\begin{xy}
,(-2,-7)*{\omega_1} ,(5,2);(5,-4)*{\bullet}**\dir{-}
?>*\dir{>}
,(5,-7)*{\omega_2}
,(-2,3);(-2,-4)*{\bullet}**\dir{-} ?>*\dir{>}
,(-2,5)*{\chi}  
,(5,2);(5,2)**\crv{(0,7)&(5,12)&(10,7)} ?>*\dir{>}
,(7,1)*{\alpha}
\end{xy}
\quad + \quad 
\begin{xy}
,(12,-7)*{\omega_2} ,(5,2);(5,-4)*{\bullet}**\dir{-}
?>*\dir{>}
,(5,-7)*{\omega_1}
,(12,3);(12,-4)*{\bullet}**\dir{-} ?>*\dir{>}
,(12,5)*{\chi}
,(5,2);(5,2)**\crv{(0,7)&(5,12)&(10,7)} ?>*\dir{>}
,(7,1)*{\alpha}
\end{xy}
\end{displaymath}
If we define 
\begin{displaymath}
[\omega_1,\omega_2] :=
\frac{\omega_1\bullet\omega_2-\omega_2\bullet\omega_1}{2}
\end{displaymath}
then
\begin{align*}
[\omega_1,\omega_2]&
=\alpha^{ij}(\partial_i\omega_{1\,k}+\partial_k\omega_{1,i})
\omega_{2,j}\mathrm dx^k+
\alpha^{ij}\omega_{1,i}(\partial_j\omega_{2\,k}+
\partial_k\omega_{2,j})\mathrm dx^k\\
&\qquad\qquad
+\partial_k\alpha^{ij}\omega_{1,i}\omega_{2,j}\mathrm dx^k=\\
&=
(\partial_k\alpha^{ij}\omega_{1,i}\omega_{2,j}+
\alpha^{ij}\partial_k\omega_{1,i}\omega_{2,j}+
\alpha^{ij}\omega_{1,i}
\partial_k\omega_{2,j})\mathrm d x^k\\
&\qquad\qquad
-\alpha^{ij}\partial_j\omega_{1\,k}
\omega_{2,i}\mathrm dx^k+
\alpha^{ij}\omega_{1,i}\partial_j\omega_{2\,k}\mathrm dx^k=\\
&=
\mathrm d\langle \alpha|\omega_1\wedge\omega_2\rangle+
{\mathcal L}_{\alpha\lrcorner\omega_1}\omega_2
-{\mathcal L}_{\alpha\lrcorner\omega_2}\omega_1,
\end{align*}
i.e., the bracket $[\omega_1,\omega_2]$ is precisely the
Magri-Koszul bracket on $1$-forms.

In particular one can recover the Jacobi identity for the Magri-Koszul 
bracket as a Ward identity (see Section~\ref{Fiorenza:trick}) by 
choosing
\begin{displaymath}
\phi = \int_0^1 \mathrm dt \int\mathrm d\theta\, 
\omega_1(x+\tilde\xi(0))\, \omega_2(x+\tilde\xi(t,\theta))\, 
\omega_3(x+\tilde\xi(1))\, \chi(\tilde\xi(\infty),\tilde\eta(\infty)).
\end{displaymath}

\subsection*{Acknowledgements}

We thank Alberto Cattaneo for having introduced us to the subject and for
useful discussion. We also thank Maciej Blaszak, Paolo
Cotta-Ramusino, Jim Stasheff and Blazej Szablikowski for their interest.

\end{document}